\documentclass[a4paper, 12pt]{article}
 \usepackage{graphicx, color}
\usepackage{graphics}
\usepackage{titling}

\pretitle{%
\begin{center}
\LARGE
\includegraphics[width=0.9\textwidth]{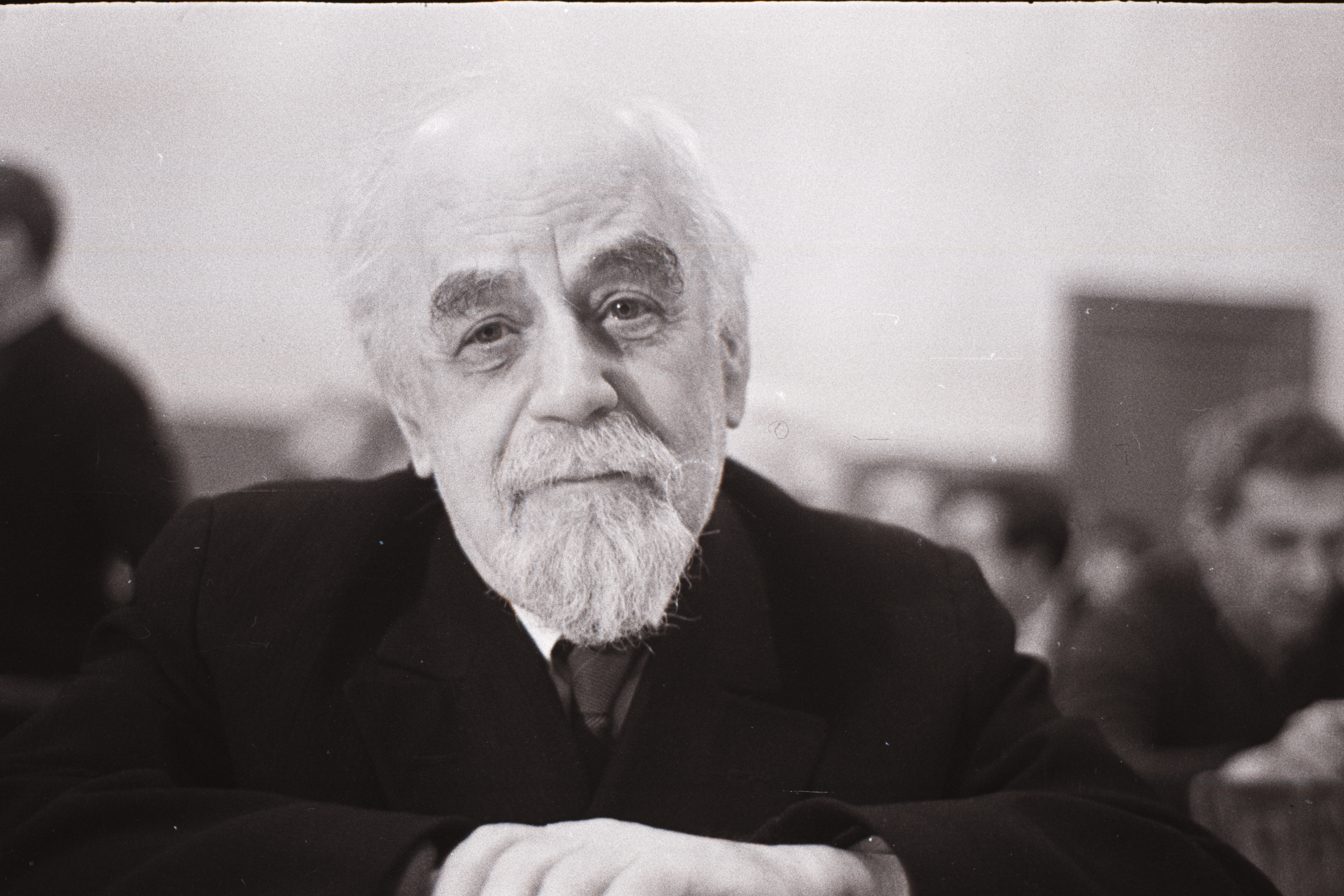}\\[\bigskipamount]
}
\posttitle{\end{center}}

\begin{document}

\title{Vladimir Ivanovich Smirnov\\
(1887-1974)}
\author{Darya Apushkinskaya and Alexander I. Nazarov\thanks{The authors are grateful to Prof. Alexander Sobolev for his valuable advice which helped us to improve the use of English in the manuscript.}}

\maketitle
The year 2017 
marked
the 130th anniversary of the prominent  Russian mathematician Vladimir Ivanovich Smirnov. 
He was a representative and 
a follower
of the famous St. Petersburg mathematical school, 
the origin of which can be traced back
to L.~Euler. 
He was  an outstanding figure of mathematical education and phenomenal organizer of science. His knowledge of physics, history, 
philosophy and music made him a real encylclopedist. 
It is hard to 
grasp the magnitude of his achievements in their entirety.
However, let's take everything in order...

Vladimir Smirnov, born on 10 June 1887 in St. Petersburg, was the youngest 
 of the ten
sons of a clergyman.
In high school, Vladimir attended the 
renowned 
Second Gymnasium, the oldest gymnasium in Russian Empire. 
Here he was luck to 
learn
mathematics from an excellent teacher Ya.V. Iodynskiy. 
Apart from 
the compulsory
classes, Iodynskiy organized a home mathematical circle where he additionally trained 
 his pupils 
for independent scientific work. Together with V.I. Smirnov, an active part in the circle took 
A.A. Friedmann\footnote{Alexander A. Friedmann (1888-1925) was a Russian physicist and mathematician best-known for his pioneering 
work on the theory of relativity. }%
\ and Ya.D. Tamarkin\footnote{Yakov (Jacob) D. Tamarkin (1888-1945) was a Russian-American mathematician who made important contribution 
to mathematical analysis and ODEs. He was a proponent and a founding co-editor of the \textit{Mathematical Reviews}. He was also a 
Vice-President of AMS in 1942-1943. }.%
\ Friendship and scientific collaboration between the circle participants continued later for many years. 

In 1905, Vladimir graduated from 
high school with 
 a 
gold medal and entered the Physics and Mathematics Faculty of St. Petersburg University. 

Smirnov's student years were marked by  anti-goverment protests and strikes in Russia. The University was 
buzzing with the political unrest.
 Lectures were held irregularly. During the strikes, the classes were cancelled completely. 
Examinations could be taken at any time 
during the 
academic 
year, by appointment. Such a system required an independent work of students with textbooks,  scientific monographs and papers, and 
 it
encouraged the creation of  circles and working groups. Among the members of the student group, founded in 
1908, were Smirnov, Friedmann, Tamarkin,  Ya.A. Shohat\footnote{Yankel (James) A. Shohat (1886-1944) was a Russian-American mathematician
 who worked on the moment problem.}%
\ and A.S. Besicovitch\footnote{Abram S. Besicovitch (1891-1970) was a Russian-British mathematician, world-known specialist in the
 function theory, winner of the Silvester Medal.}.\ These students studied almost without  guidance from the university. However, some
 professors observed their development with interest.

In 1910, Smirnov graduated from the University with  a 
First Class Diplo- ma.
During the next two years, he taught mathematics 
 at
one of the best private schools in St. Petersburg. Also at this time, Smirnov and his friends organized 
an
informal seminar, where 
they taught themselves 
mini-courses in various 
fields
of mathematics and mechanics. These mini-courses included, in partucular, 
theory of surfaces, complex analysis, potential theory, spherical functions.

In 1912, 
Smirnov was 
admitted to the University
as a graduate student  
of Prof. V.A. Steklov\footnote{Vladimir A. Steklov (1864-1926) was a prominent Russian mathematician, mechanician and physicist. 
Nowadays, the Mathematical Institutes in Moscow and St. Petersburg are named after Steklov.}.\ 
The study program was quite intensive. 
As a master student\footnote{The Master's Degree in Russian Empire corresponds to today's PhD.}, Smirnov had to sit four exams:
differential equations and mathematical physics, 
number theory, 
probability theory and 
theoretical mechanics. 
The 
reading
list 
included the works of Fourier, Sturm, Liouville, Floquet, Poincar\'{e}, Lyapunov and Steklov. 

In addition to the 
recommended reading,
Smirnov was interested in the latest results on analysis such as 
theory of integral equations.
Steklov was skeptical about these new directions, so Smirnov  studied  the works of Frechet and the Riesz-Fisher Theorem on his own.

In 1913 Smirnov published (jointly with Friedmann) his first paper on the oscillatory discharge of a capacitor.
The same year Vladimir Ivanovich married Ekaterina Nikolaevna Gorbunova, who 
was a gymnasium history teacher.
Also at that time  Smirnov began his teaching activities at universities. From 1912 he taught at  St. Petersburg Mining Institute 
and at  St. Petersburg Institute of Railway 
Engineering, 
and from 1916 he worked as Prof. A.V. Vasil'ev's\footnote{Alexander V. Vasil'ev (1853-1929) was a Russian mathematician, known for 
his activity as historian of mathematics and organizer of scientific life.} 
assistant
at  St. Petersburg University. 

So, in the beginning Smirnov's life ran quite smoothly: high school, interest in mathematics, university, 
work 
under the guidance of the well-known mathematician Steklov, marriage, promotion.
The future was promising to be 
happy and secure.
However, the First World War (1914-1918), the Revolution (1917) and  the Civil War (1918-1920) 
have completely destroyed these hopes.
There was not enough food, heating...

In 1918, Smirnov 
defended his 
thesis \textit{``The  inversion problem of a linear second-order differential 
equation with four singular points''}. 
He had chosen the topic of his 
Thesis by himself. 
Steklov did not find
this topic 
very attractive, but wrote a brilliant 
report
and emphasized that Smirnov's work had opened new directions
in the analytic theory of differential equations.

During the difficult post-revolutionary period many scientitsts  were forced to leave hungry Petrograd\footnote{In 1914, the city 
name was changed from St. Petersburg to Petrograd. In 1924 it was renamed to Leningrad, and in 1991 back to St. Petersburg.}. 
In this situation, V.A.~Steklov 
helped 
his 
students 
to find positions 
at
provincial universities. Thus, in the autumn of 1918
Vladimir Ivanovich got a job as a 
privatdozent 
at the newly founded Tavria University in Simferopol, Crimea. By this time  his family, i.e., his wife and her two daughters from the first marriage, also lived in Crimea, in Balaklava. 

Life in 
Crimea 
turned out to be
even scarier as in Petrograd. 
By 1920 
the regime in Simferopol has changed seven times, 
and ordinary people 
kept disappearing.
The greatest bloodshed  
occurred
in 1920 when the  Red Army had captured 
Crimea\footnote{More than 50 000 people were shot in Crimea by the end of 1920.}. The winners regarded the university with liberally oriented scientists as an enemy lair.
Smirnov, like many others, received a summon to the so-called ``court''. Fortunately, a ``judge'' allowed him to go home.
At the same time
the wife of Vladimir Ivanovich was 
executed
by a court ruling\footnote{Smirnov supported the daughters of Ekaterina Nikolaevna  for many years.}.

Smirnov 
recounted this tragedy 
to Steklov and also informed him about  
 the catastrophic situation 
at the 
 Tavria University, 
 asking for help with moving back to Petrograd.
This long-awaited trip   
took place in 1921. Smirnov returned to his
\textit{alma mater} and never left the University again. 
At this
moment,  the long-term activity of Vladimir Ivanovich  
as an organizer of mathematical life in Petrograd-Leningrad began. 

Smirnov developed a fundamentally new mathematical course  for  physicists,  
encompassing 
the modern knowledge and methods.
He gave lectures for physics students that 
spanned an unheard-of range of mathematical disciplines.
In 1925, Smirnov founded the department 
of Theory of Functions of a Complex Variable. 
Also during these years, Smirnov 
established
several 
regular
scientific seminars: on the theory of functions of a complex variable, 
analytic
theory of differential equations, 
approximation
methods, 
functional analysis. These seminars attracted promising enthusiastic 
young Leningrad mathematicians, and thus they paved the way to the creation
of world-known scientific schools in all these areas.
 
Vladimir Ivanovich was in touch with all major achievements in science. In 1920's he was one of the organizers and active members of the Leningrad Mathematical Society, which unfortunately 
was dissolved in 1930 due to political reasons. Smirnov 
was also at the heart of the Second All-Union Congress of Mathematicians (1932) in Leningrad, and we are also indebted to him for 
making the publication of the  proceedings of the congress possible.

From 1929 to 1935 Smirnov was the head of
the theoretical department of the 
Seismology 
Institute. Here, together with S.L.~Sobolev\footnote{Sergei L. Sobolev (1908-1989) was an outstanding Soviet mathematician working 
in mathematical analysis and PDEs. He introduced the notions (Sobolev spaces, distributions, etc.) that are now fundamental for 
several areas of mathematics. Nowadays, Institute of Mathematics in Novosibirsk is named after Sobolev.}, he applied a  method of functionally invariant solutions for solving some complicated problems related to 
wave propagation.

In 1931, Smirnov 
became the deputy director
of the Institute of Mathematics and Mechanics, organized 
on his 
initiative  
at the University. Two years later, he founded
the department 
of Higher Mathematics at the Faculty of 
Physics. Smirnov  
remained its head
for more than 40 years. 

1934 was a special year for Vladimir Ivanovich: he married Elena Prokopi- evna Ochlopkova, 
who was his graduate student. After marriage she did not continue her research activity but worked as an Assistant Professor at the University. In 1935 
the
 Smirnovs got their son Nikita.

The 
peaceful
life was 
destroyed
by the World War II. 
At the end of August 1941 Smirnov and his family, together with 
other University staff
were evacuated  from Leningrad to Yelabuga\footnote{Yelabuga is a small town in the Republic of Tatarstan, located 200 kilometers east from Kazan.}. 
Conditions of daily life here were very hard. 
Most of the  families had to content themselves with just one room to live.
Also, the University employees 
had to provide themselves with 
firewood
and to grow vegetables 
in a small garden 
to save themselves from starvation.

Despite all the hardships,  
Smirnov continued to work actively. He taught mathematics at the Pedagogical Institute, 
and served as a high scholl examiner.
In addition, Vladimir Ivanovich organized an 
aerodynamics group, 
which had carried out a series of 
studies
on ballistics under his supervision. In 1943, Smirnov was elected a Member of the Academy of Sciences of the USSR.

Upon his
return to Leningrad in 1944, Smirnov found the Faculty of Mathematics 
 and Mechanics of Leningrad University in a 
pitiful  
 state. The 
 numbers of students were very low, professors even lower.
For this reason, during the next 12 years Smirnov 
was in charge 
(successively
or sometimes even  simultaneously)  of several  
departments:
Elasticity Theory, 
Hydroaerodynamics, Theoretical Mechanics, Theory of Functions of a Complex Variable, and  
Mathematical Analysis.
As soon as a 
 suitably qualified 
candidate appeared, Vladimir Ivanovich would step down and 
hand over the responsibility  
to the young colleague.

Being for many years the only mathematician member of the Academy of Sciences 
 in Leningrad, Smirnov had to communicate 
 papers 
 in  all fields of mathematics to the \textit{Doklady of Academy of Sciences}. 
 He never treated this duty formally, and if necessary, he  helped  to 
improve  the papers. There were cases when 
papers were completely rewritten with Smirnov's help.

Already in  the early 1920's, Smirnov 
and Tamarkin began to write a textbook for technicians and physicists. The first two volumes were published in 1924 and 1926, respectively. However, the joint work did not continue for long because Tamarkin soon 
fled from the  country.
For further editions Smirnov extended and revised the material alone. 
Since 1930, ``\textit{A Course of Higher Mathematics}''  began to appear only under the name of Smirnov.

By 1947, the \textit{Course} has become a 
veritable
mathematical encyclopedia in five volumes. 
It presents
a systematic treatment of Calculus, Algebra, Differential Geometry, Ordinary Differential Equations, Vector Analysis, Fourier Series,
 Complex Analysis (including functions of several complex variables), Special Functions, Theory of Group Representations, Calculus of
 Variations, Integral Equations, Partial Differential Equations and Boundary Value Problems. The 5th volume was the first Russian
textbook on Functional Analysis.  In 1948, the \textit{Course} was awarded the Stalin Prize. 

Smirnov continued to refine and 
supplement
the \textit{Course} throughout his life. Actually this activity 
stretched
for over 50 years. It was truly Smirnov's life's work!
The\textit{ Course}  was reissued many times (in particular,  up to now the first two volumes 
were reprinted
24 times) and translated into 8 languages.  
Nowadays, these 
books can be found in almost every university library
around the world.

In the fall of 1947, Vladimir Ivanovich
initiated
 a scientific seminar on mathematical physics 
bringing together  
 a major part of the Leningrad mathematicians working in the field of partial differential equations and their applications. 
Almost all Soviet experts in PDEs  viewed
a possibility 
to give a talk  at this seminar as an honour. 
Among the speakers of the seminar there were also many famous foreign mathematicians such as R. Courant, J. Leray, P. Lax, K.O. Friedrichs etc. Now this seminar is named after V.I. Smirnov. In 2017, it celebrated 
its 70th anniversary.
 
In 1956, on 
Smirnov's
initiative 
the 
department
of Mathematical Physics was 
established
at the Faculty of Mathematics and Mechanics. 
The 
purpose
of this department
was the training of highly qualified specialists in 
modern problems of PDEs. 
Smirnov hired for
teaching both well-known scientists (e.g., S.G. Mikhlin\footnote{Solomon G. Mikhlin (1908-1990) was a famous Soviet mathematician  working in analysis, integral equations and computational mathematics. He is best known for the introduction of the concept of ``symbol of a singular integral operator''.}, 
O.A. Ladyzhenskaya\footnote{Olga A. Ladyzhenskaya (1922-2004) was a prominent Soviet and Russian mathematician known for her work on PDEs (especially Hilbert's 19th problem) and fluid dynamics. President of the St. Petersburg Mathematical Society (1990-1998).}) and young talented mathematicians. Smirnov was 
the
head of this 
department
and of the 
Higher Mathematics department at the Physics Faculty) 
until
the end of his life. 

In the mid-fifties, after Stalin's death, the idea  of the revival of  Leningrad Mathematical Society 
was again in the air.
In 1957 Smirnov 
established
the Leningrad General Mathematical Seminar. Two years later,
on the basis of this seminar  the Mathematical Society was restored, and Vladimir Ivanovich was elected its  Honorary President. 

Smirnov kept working actively until the age 85. In 1972 
he was struck down by a sudden illness which took his life in two years.
Vladimir Ivanovich passed away on 11 January 1974. He was buried 
on 
the  Komarovo Cemetery\footnote{Komarovo is a small summer-resort in 45 km away from Leningrad. Smirnov, among some other members of the Soviet Academy of Sciences, had summer residence (``dacha'') there.}. 
$$
***
$$

Smirnov's scientific results 
 are 
  overshadowed
 by his tremendous  educational  and organizational activities. 
Nevertheless, one should not forget that 
 he is the 
 author of several remarkable results in various fields of analysis. 

The first important cycle of Smirnov's works (including his 
Thesis) is connected with the analytical theory of  ordinary differential equations. 
Among  these works we mention only the 
ones that are concerned
with the problem of inversion for a linear second-order differential equation with four singular points. For the case of three singular points, such a problem was comprehensively studied by K.H.A. Schwarz. Smirnov gave a complete solution of the problem with four singularities.
His investigations were carried out in a quite elegant analytical way.

A large number of Smirnov's papers are devoted to the study of the boundary properties of analytic functions and theory of approximations. He found the canonical factorization of a function from Nevanlinna class\footnote{Later this factorization was 
rediscovered
by R.H. Nevanlinna. Sometimes it is called Nevanlinna factorization.}. 
He also singled out special subclasses of the Nevanlinna class. 
Now these subclasses are called Smirnov's classes. They are the best known and most thoroughly 
studied generalizations of the Hardy classes. 
Smirnov's classes emerge  in natural way when tackling various analytic questions, including those from 
approximation theory.
When studying approximations,
Smirnov introduced a class of domains 
that
later became known as Smirnov domains.

In 
the
early 1960's, 
jointly with N.A. Lebedev\footnote{Nikolai A. Lebedev (1919-1982) was a Soviet mathematician who worked on the theory of  functions of a complex variable.}, Smirnov wrote 
the 
 monograph  
``\textit{Functions of a complex variable: constructive theory}'' which 
contained modern results on the subject, including the authors' contribution.
The book was translated into English and highly appreciated by the experts.

The third essential group of Smirnov's results relates to the propagation of sound and elastic waves. He discovered the existence of the so-called functional invariant solutions of the wave equation (families of solutions depending on functional parameters). As a result, all physically relevant
fundamental solutions of the wave equation with three space variables were constructed. Singular solutions of the three dimensional dynamical system of elasticity theory possessing axial symmetry were constructed as well.
The works of Smirnov and 
Sobolev brought the Soviet school of theoretical seismology  to the 
forefront of this research in the world. 

Considering the  mixed boundary value problems for the wave equation and for the dynamical system of  elasticity theory in the ball,
 Smirnov suggested a method which was later called the Smirnov method of  incomplete separation of variables\footnote{also known as  the space-time triangle diagram (STTD) technique.}.

$$
***
$$

Smirnov was a brilliant university teacher and lecturer. His lectures were not only elegant and comprehensible. In his exposition, the functions, integrals, and equations behaved as animate creatures.

His 
lectures on general courses 
Smirnov began with classical examples, then explained the basic notions and ideas.
After this, he 
would change the tone of his voice, and in the last ten minutes of the lecture 
would give 
two or even three times more material as before, with the possible generalizations 
of both 
mathematical and 
physical nature.
Each student could decide for himself how far to follow the lecturer. 
Note that
Smirnov's lectures were also attended by 
people
who were not formally university students.
\begin{figure}[!htbp]
\centering
\includegraphics[width=0.9\textwidth]{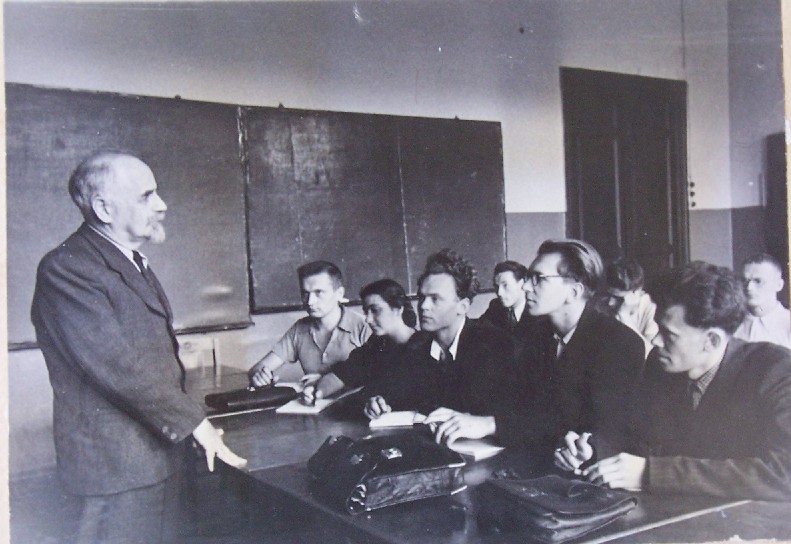}
\caption{V.I. Smirnov during the lecture (1951)}
\end{figure}

According to \textit{Mathematics Genealogy Project}, Vladimir Ivanovich had only 12 PhD students and Postdocs: N.P. Erugin\footnote{Nikolay
 P. Erugin (1907-1990) was a Soviet mathematician working in the field of ODEs.}, G.M. Goluzin\footnote{Gennadiy M. Goluzin (1906-1952) 
was a Soviet mathematician known for his works in geometric theory of functions of a complex variable.}, V.P. Havin\footnote{Victor 
P. Havin (1933-2015) was a prominent Soviet and Russian mathematician known for his works in analysis. Among his students there are 
many famous analysts including the Fields medalist Stanislav K. Smirnov (2010).}, L.V. Kantorovich\footnote{Leonid V. Kantorovich
 (1912-1986) was an outstanding Soviet mathematician working in the fields of functional analysis, computational mathematics and
 mathematical economy. He introduced the notion of normed vector lattices, which are called \textit{K-spaces} in his honor. He was 
the winner of  the Nobel Memorial Prize in Economics (1975) for ``contribution to the theory of optimum allocation of resources''.}, 
P.P. Korovkin\footnote{Pavel P. Korovkin (1913-1985) was a Soviet mathematician working in the field of analysis.}, 
V.D. Kupradze\footnote{Viktor D. Kupradze (1903-1985) was a Soviet mathematician known for his works in the elasticity theory 
and integral equations.}, I.A. Lappo-Danilevsky\footnote{Ivan A. Lappo-Danilevsky (1896-1931) was a Soviet mathematician known 
for his fundamental results on  linear differential equations and analytic functions of matrices.}, S.M. Lozinsky\footnote{Sergey 
M. Lozinsky (1914-1985) was a Soviet mathematician known for his fundamental results on ODE s and computational mathematics. 
President of the Leningrad Mathematical Society (1965-1985).}, S.G. Mikhlin,
A.L. Shahinyan\footnote{Artashes L. Shahinyan (1906-1978) was a Soviet mathematician known for his works in the theory of functions 
of a complex variable.}, S.L. Sobolev and V.A. Yakubovich\footnote{Vladimir A. Yakubovich (1926-2012) was a prominent Soviet and 
Russian mathematician known for his works in control theory, winner of the IEEE Control Systems Award (1996).}. 
However, this number was immeasurably smaller than the number of people who 
profited from
Smirnov's supervision, consultations, and advice.

In
Smirnov's own words,  the main goal of his 
life was to carry 
the knowledge, received from his great teachers, through 
the years of wars and ruin, and 
pass 
it on to new generations.

Smirnov's extensive knowledge  is well illustrated 
by the following example
: in 1950's he looked through \textbf{\underline{all}} the PhD theses 
in
mathematics and mechanics, defended in Leningrad.
A deep 
general
understanding of mathematics as 
one entity, explains why Smirnov supported and promoted even those mathematical areas 
where he had never worked himself.
In particular, such 
areas include functional analysis and spectral theory of operators.

Mathematical and general 
Smirnov's encyclopedic knowledge
and his proficiency in foreign languages made him a key person in such an important matter as 
new acquisitions for
the Library of the Academy of Sciences. 
For many years, Vladimir Ivanovich was a consultant 
on mathematical and physical literature for the acquisition department of the library.

$$
***
$$
Vladimir Ivanovich was a rare example of a scientist who combined his own research  with  
professional activity as historian of his science. 

Smirnov worked actively in the Commission on the History of Physical and Mathematical Sciences of the Academy of Sciences\footnote{In 1953,
 this Commission has become a part of the Institute of the History of Natural Science and Technology.}, and from 1951 to 1953 he served as
 its chairman.
He was   a member of the Editorial Boards and the author of several survey articles in the multivolume editions of 
\textit{``The Complete Works by P.L. Chebyshev\footnote{Pafnuty L. Chebyshev (1821-1894) was an outstanding Russian mathematician and mechanician, founder of the St. Petersburg mathematical school,  who is remembered primarily for his works on number theory, probability
 and approximation theory.}''} (1946-1951), 
 \textit{``The History of the  Academy of Sciences of the USSR''} (1958-1964), and \textit{``The History of Homeland Mathematics''}
 (1966-1970). Vladimir Ivanovich was also the editor of the first Russian edition of \textit{``The Selected Works by M.V. 
Ostrogradsky\footnote{Mikhail V. Ostrogradsky (1801-1862) was an outstanding Russian mathematician, mechanician and organizer of mathematical
 education, who made the essential contributions in many fields of mathematics and physics.}''} (1958), and of the fundamental collective volume \textit{``Mathematics at the University of  Petersburg-Leningrad''} (1970), for which Smirnov wrote several chapters too.

Vladimir Ivanovich 
 had invested a tremendous effort into reviewing and publishing
manuscripts 
of
outstanding Russian mathematicians. Among them we mention D. Bernoulli\footnote{Daniel Bernoulli (1700-1782) was an outstanding Swiss
 mathematician and physicist, one of the founders of mathematical physics. He worked in St. Petersburg from 1725 till 1733. In this period,
 Bernoulli prepared his world known monograph \textit{``Hydrodynamica''}. }, A.M. Lyapunov\footnote{Alexander M. Lyapunov (1857-1918) was an
 outstanding Russian mathematician and mechanician. He created the stability theory of the dynamical systems and established several
 breakthrough results in mathematical physics and probability. }, A.N. Krylov\footnote{Alexey N. Krylov (1863-1945) was an outstanding
 Russian and Soviet naval engineer and applied mathematician, a founder of Russian school of theoretical naval architecture.},
I.A. Lappo-Danilevsky. 

The last case is, probably, unique in the history of science.  A former student of Smirnov, Lappo-Danilevsky passed away at the age 34.  Within 4 years after that,\  Vladimir Ivanovich\  (together\  with \ N.E. Kochin\footnote{Nikolai E. Kochin (1901-1944) was a famous Soviet
 mathematician, one of the founders of theoretical meteorology.})\  studied \ his manuscripts and rough drafts, filled in all the gaps, and 
eventually
published 12 papers by Lappo-Danilevsky. This work can be likened to
that
of a fine art restorer.

A very special place in Smirnov's historical investigations
 is reserved for
 the study of the scientific heritage of  L. Euler\footnote{Leonhard Euler (1707-1783) was a great Swiss and Russian mathematician and
 physicist, who made fundamental discoveries in all branches of mathematics as well as in various natural sciences.}. At the end of the 1950's Vladimir Ivanovich 
initiated 
 a systematic study of the extensive St. Petersburg archives of Euler including  hand-written 
 notes 
 and correspondence. In particular, to study Euler's notes on 
 number theory,  Smirnov 
 appointed his own PhD student G.P. Matvievskaya\footnote{Galina P. Matvievskaya (born 1930) is a Soviet and Russian expert in history of mathematics and oriental studies.}. 

As a result of these massive investigations, several volumes of Euler's works, including a 438 pages annotated index,  were published. Smirnov also served as a president of Soviet part of joint Euler's Committee organized by  Schweizerische Naturforschende 
Gesellschaft\footnote{since 1988 Swiss Academy of Natural Sciences.} and Academy of Sciences of the USSR.

Vladimir Ivanovich  
also penned
a large number of 
biographical 
essays on Russian mathematicians such as Chebyshev, Steklov, Friedmann, 
N.M. G{\"u}nter\footnote{Nikolai M. G{\"u}nter (1871-1941) was a famous Russian and Soviet mathematician known for his works in potential theory and PDEs.} and many others.
These essays imprinted the living images of remarkable people.

For many years Smirnov 
served as the 
President of the Scientific Council of the Academic Archives 
and a member of the Scientific Council of Academic Library. 
He was the author and editor of several extensive bibliographic indices. 
His colleagues pointed out that
 a task requiring the work of a research group for its implementation was often carried out by Smirnov alone and was done very well and within the shortest time.
 
Smirnov's  
 deepest knowledge  of 
 general history and 
of  
 the history of science, as well as his special care in analysing  factual material, 
 amazed the contemporaries. Smirnov's articles on topics related to the history of science were so accurate in all the details and facts that they can be used 
 as reliable sources for necessary references.

$$
***
$$

The most important hobby of Smirnov throughout all his life was music. He knew music professionally, played 
piano well, and had a good musical memory. He often took to the concerts, that he attended, the scores and compared the  played music with the original. 

Every two weeks Vladimir Ivanovich  
held
excellent musical recitals at home, during which he often played popular symphonies (for example by G. Mahler) adapted for piano 
four-hands. For a long time, his music partner was D.K. Faddeev\footnote{Dmitry K. Faddeev (1907-1989) was a prominent Soviet 
mathematician known for his works in algebra, number theory and computational mathematics.}. When Smirnov and Faddeev  played four-hands,
 no one was allowed to speak.
Smirnov was  
immensely
proud that at his home the listeners could enjoy 3 symphonies and 2 quartets in one evening.
Vladimir Ivanovich also played duets with D.D. Shostakovich\footnote{Dmitry D. Shostakovich (1906-1975) was a Soviet composer and pianist,
 regarded as one of the greatest composers of the 20th century.}, 
who at the turn of the 1950's spent  a summer 
in Komarovo at a dacha located next to Smirnov's dacha. 

Only in his last years  Smirnov was forced to restrict himself to listening to musical records.

Of literature, Smirnov liked the classical one. His favorite writer and the one he found most congenial was Dostoevsky.
In the last decades of his life, Vladimir Ivanovich preferred to read  memoirs. Also, Smirnov read a lot of philosophical literature 
(in particular, Spinosa) in the original languages.

Vladimir Ivanovich had a great love for Russian architectural heritage, expecially the churches and monasteries of the North-West part of Russia. Almost every summer he visited some secluded corner, favoured by the masters of past centuries.

Smirnov liked brisk
long  distance  walks.  They helped him to 
remain cheerful and be efficient
, and also served as a kind of healing  for his not quite healthy heart.

Another hobby of Vladimir Ivanovich was playing cards, namely  ``vint''\footnote{Vint (literally ``screw'') is a complicated card game of the bridge class.}. One of his regular card partner was N.M. G{\"u}nter.

$$
***
$$

Vladimir Ivanovich was a deeply religious person, an Orthodox Christian.  
During the Soviet period, 
practising religion 
was very difficult, and sometimes even dangerous. 
Nevertheless, he regularly 
made donations to the Church,
and, for a long period,  was a member of ``The Twenties'' 
(Parish Council) of the Prince St. Vladimir's Cathedral in Leningrad. 
At the same time, Smirnov never stressed his 
religiousness in every day contacts with people.

In the years of Stalin's terror, Smirnov was not afraid to  write letters to the authorities (prosecutors, ministers, politicians) and to speak in defense of arrested colleagues.  After one of these speeches, 
an NKVD\footnote{NKVD (the People's Commissariat for Internal Affairs) was the leading Soviet secret police organization from 1934 to 1946.} officer paused ominously and then said to him: 'I see, you are 
quite a brave man, Vladimir Ivanovich!' In addition, Smirnov supported the families of 
prosecuted
people both financially and morally.

Also, many persons were obliged to Smirnov 
for the opportunity to remain mathematicians.
 If a talented student was not accepted for postgraduate study 
for
  political reasons or a mathematician was sent to work, where he could not use his professional knowledge, Vladimir Ivanovich fought for such people like a lion. 
 
The humanism, 
integrity, 
incorruptibility and fortitude were Smirnov's outstanding features, determining his influence.
Smirnov was also  admired by all who knew him for his attractive personal qualities, his rare charm and exceptional modesty. His social circle included not only scientists, but also famous writers, artists, composers and musicians.

During his long life, this extraordinary man performed many good, wise and useful deeds. Vladimir Ivanovich lived and worked trying to give as much and take as little as possible, and this influenced people around him stronger than any orders or sermons on morality.

\end{document}